\documentclass[a4paper,times,3p]{elsarticle}
\pdfoutput=1
\usepackage[utf8]{inputenc}

\usepackage[hidelinks]{hyperref}

\usepackage{booktabs}

\usepackage[font=footnotesize,labelfont=bf]{caption}

\usepackage{todonotes}
\usepackage{algpseudocode}
\usepackage{setspace}

\usepackage{microtype}

\usepackage{aicescover}
\usepackage{amsmath}
\usepackage{amssymb}
\usepackage{marvosym}
\usepackage{gnuplot-lua-tikz}
\usetikzlibrary{positioning,shapes,arrows,backgrounds,calc,fit}
\usepackage{listings}

\definecolor{mygreen}{rgb}{0,0.6,0}
\definecolor{mygray}{rgb}{0.5,0.5,0.5}
\definecolor{mymauve}{rgb}{0.58,0,0.82}

\definecolor{RWTHBlueDark}{RGB}{0,84,159}
\definecolor{RWTHBlueLight}{RGB}{142,186,229}
\definecolor{RWTHRed}{RGB}{204,7,30}
\definecolor{RWTHGreen}{RGB}{87,171,39}
\definecolor{RWTHGreen2}{RGB}{189,205,0}

\newcommand{\D}[2]{\frac{\partial #1}{\partial #2}}

\renewcommand{\vec}[1]{\mathbf{#1}}
\newcommand{\mat}[1]{\mathbf{#1}}
\newcommand{\bm}[1]{\mathbf{#1}}

\newfont{\logobold}{logobf10 scaled\magstep2}
\newcommand{\assembly}{\mathop{\mbox{\logobold A}}}

\aicescovertitle{Automatic implementation of material laws: Jacobian calculation in a finite element code with TAPENADE}
\aicescoverauthor{Florian Zwicke \and Philipp Knechtges \and Marek Behr \and Stefanie Elgeti}

\def\statement{\begin{minipage}[t]{.75\textwidth}
       NOTICE: This is the author's version of a work that was published in
Computers \& Mathematics with Applications 72 (2016) 2808--2822,
\href{http://dx.doi.org/10.1016/j.camwa.2016.10.010}{DOI:10.1016/j.camwa.2016.10.010}.
\\

\copyright \, 2016. This manuscript version is made available under the
CC-BY-NC-ND 4.0 license \\
\url{http://creativecommons.org/licenses/by-nc-nd/4.0/}.
       \end{minipage}}

\makeatletter
\def\ps@pprintTitle{%
     \let\@oddhead\@empty
     \let\@evenhead\@empty
     \def\@oddfoot{\footnotesize\itshape
       \statement\hfill\today}%
     \let\@evenfoot\@oddfoot}
\makeatother

\journal{}

\begin{document}
\aicescoverpage
\begin{frontmatter}

\title{Automatic implementation of material laws: Jacobian calculation in a finite element code with TAPENADE}

\author{Florian Zwicke\corref{mycorrespondingauthor}}
\cortext[mycorrespondingauthor]{Corresponding author}
\ead{zwicke@cats.rwth-aachen.de}
\author{Philipp Knechtges\corref{}}
\ead{knechtges@cats.rwth-aachen.de}
\author{Marek Behr\corref{}}
\ead{behr@cats.rwth-aachen.de}
\author{Stefanie Elgeti\corref{}}
\ead{elgeti@cats.rwth-aachen.de}
\address{Chair for Computational Analysis of Technical Systems (CATS)\\
Center for Computational Engineering Science (CCES), RWTH Aachen University, 52056 Aachen, Germany}

\begin{abstract}
In an effort to increase the versatility of finite element codes, we explore the possibility
of automatically creating the Jacobian matrix necessary for the gradient-based solution of
nonlinear systems of equations.
Particularly, we aim to assess the feasibility of employing the automatic differentiation
tool TAPENADE for this purpose on a large Fortran codebase that is the result of many
years of continuous development.
As a starting point we will describe the special structure of finite element codes and
the implications that this code design carries for an efficient calculation of the
Jacobian matrix. We will also propose a first approach towards improving the efficiency of
such a method.
Finally, we will present a functioning method for the automatic implementation of the
Jacobian calculation in a finite element software, but will also point out important
shortcomings that will have to be addressed in the future.
\end{abstract}


\end{frontmatter}

%\linenumbers

\section{Introduction}

In the field of numerical simulations, a central task
that needs to be accomplished is the solution of
partial differential equations (PDEs).
When such equations are discretized, the result is a system
of equations. For a set of nonlinear partial differential
equations, the discretized equations are also nonlinear.

A very popular method for the solution of nonlinear systems
of equations is the Newton--Raphson iteration scheme \cite{Ypma1995}.
If a nonlinear system of equations for a vector of unknowns $\vec{u}$ is given by

\begin{equation}
  \vec{F}(\vec{u}) \; = \; \bm{0} \, ,
\end{equation}

then, starting from an initial guess $\vec{u}_0$, one iteration
step of the Newton--Raphson scheme is given by

\begin{equation}
  \D{\vec{F}(\vec{u}^k)}{\vec{u}} \cdot \Delta \vec{u}^{k+1} \; = \; - \vec{F}(\vec{u}^k) \, ,
\end{equation}

which defines a linear system of equations.

There are many different methods for the solution of the resulting linear equation
systems. Some of these require the full Jacobian matrix
$\D{\vec{F}(\vec{u}^k)}{\vec{u}}$. Others, especially Krylov-based methods
such as GMRES \cite[\S 3]{Saad1986}, require matrix-vector products,
which would be directional derivatives of $\vec{F}(\vec{u})$ in this context.
These methods can be implemented in such a way that the Jacobian matrix does
not need to be calculated explicitly. They are called matrix-free methods
\cite[\S 3]{Hovland2001}, \cite{Brown2010}. However, the combination
of such a method with a general-purpose preconditioner usually requires the full matrix again.
For this reason, we will only consider the accumulation of the full Jacobian
matrix.

It is very common to implement the calculation of this matrix directly
into the program code. In these cases, the residual $\vec{F}(\vec{u})$ first
needs to be differentiated with respect to the vector of unknowns $\vec{u}$,
which could, for instance, be done manually. The process of differentiating
the residual and subsequently implementing the derivative depends entirely
on the differential equations in use and needs to be repeated whenever these
equations are changed.
Using this method, the task of trying out new simulation models may require
considerable effort. Additionally, a manual derivation and implementation
of the derivative invites mistakes that may be very hard to find
later on. As an alternative, there is also the option of obtaining the derivative
by an automated procedure.

There exist different approaches of this kind,
and they typically
require a program function that calculates the residual $\vec{F}(\vec{u})$
to work with.
The oldest and simplest method would be a finite-difference (FD) approximation.
This just requires repeated calls to the residual function with a perturbed input
$\vec{u} + \Delta \vec{u}$ to obtain a directional derivative. This
method has been used for the estimation of Jacobian matrices for a long time,
and some techniques for increased efficiency have also been developed. For example,
if the sparsity pattern of the Jacobian is known beforehand, an optimized
seed matrix can be used to calculate a compressed Jacobian \cite{Curtis1974}.

In the last decades, there has been a lot of interest in using the
method of automatic differentiation (AD) \cite{Kedem1980,Griewank1989}
for this purpose. Recent examples of such approaches are described in \cite{Bucker2010}
and \cite{Reynolds2012}.
The methods for sparsity exploitation that were
developed for finite-difference approaches can easily be transferred to the similar
forward mode of automatic differentiation, but more complex possibilities have also
been researched \cite{Griewank1991,Griewank1990,Coleman1998,Naumann2008}.

The most important advantage of AD over FD is that the  resulting derivatives obtained this
way are generally more accurate, since they do not suffer from possible cancellation errors
that limit the number of significant digits.
Unfortunately, this improved accuracy comes at the
cost of a much more difficult and laborious implementation.
Our aim in this paper is to investigate the suitability and efficiency of using AD
for the Jacobian calculation specifically in the context of finite element (FE) codes.

We decided to use the forward mode of automatic differentiation, since we
expect this to be more efficient for the calculation of a square
Jacobian matrix (see Section~\ref{sec:ad}).
Our research focuses on the implementation of the Jacobian calculation
in a specific finite element code written in Fortran. Of the different
techniques that are available for the implementation of automatic differentiation,
we consider source code transformation to be more suitable for Fortran codes
(see Section~\ref{sec:ad}).

There are a number of programs available that support automatic differentiation
by code source transformation in Fortran codes. One of the earliest was ADIFOR
\cite{Bischof1992,Bischof1992b}, which seems
to have been under development mostly in the 1990s. Another early implementation
was Odyss\'ee \cite{Rostaing1993}, which was later replaced by TAPENADE
\cite{TapenadeRef13}. Other software includes OpenAD/F \cite{Utke2008}.
The present work will focus on the investigation of TAPENADE as a tool for the
automatic differentiation.
This decision is motivated primarily by the long
history of this software and its widespread use in the scientific community.

One should note that even though we made some specific choices,
many of the findings we present in this document can be
applied universally, i.e., they apply for arbitrary modes or techniques of
automatic differentiation or even for arbitrary software.

Our principal aim of calculating a Jacobian matrix via forward mode AD
can easily be achieved using repetitive calculations of directional derivatives.
This is analogous to the approach used with finite-difference methods.
In cases where the Jacobian matrix is sparse, more efficient techniques can
be applied, provided that the sparsity pattern is known beforehand (see,
e.g., \cite{Curtis1974}, where this is discussed in relation to the finite differences
method).
We will show that the Jacobian matrix resulting from the finite element
discretization is indeed sparse (see Section~\ref{sec:s3implprops}).
Since its sparsity pattern depends on the
computational mesh, conventional methods of sparsity exploitation are not
suitable. We will describe an alternative way of calculating the matrix
efficiently, that utilizes structures that are commonly found in finite element
codes.

We will provide some background on relevant aspects of
automatic differentiation in Section~\ref{sec:ad}.
Section~\ref{sec:software} gives an overview of the software programs in use
as well as the difficulties arising from their combination.
We will describe the development of a suitable procedure for the automatic creation
of the Jacobian matrix in a finite element code in Section~\ref{sec:s3jacfe}, and present
the results that we obtained when using it on a viscoelastic flow problem in Section~\ref{sec:example}.

\section{Automatic differentiation with TAPENADE}

\subsection{Automatic differentiation}
\label{sec:ad}

The term automatic differentiation describes the differentiation of program functions,
when it is carried out by computer algorithms. It is based on the idea that certain
program functions can be viewed as sequences of mathematical expressions \cite{Kedem1980}.
Assuming existence, derivatives of these individual
expressions can be calculated and used in the chain rule to determine the full function derivative
 \cite{Griewank1989}.

Concerning the order in which the individual derivatives are calculated,
one generally distinguishes between a forward and a reverse
mode \cite{Griewank1989,Barthelemy1995}.
In forward mode, the calculation of the derivatives happens in the same
order as the execution of the original statements. This means that at an
arbitrary point during the calculation, the derivatives of all relevant
intermediate variables
with respect to the input variable are known. In reverse mode, on the other
hand, the operations are carried out in reverse order. Therefore, the quantities
that are available during the calculation are the derivatives of the output variable
with respect to the relevant intermediate variables.
In analogy to the matrix multiplication, where the transposition also inverts the order
in which the matrices are multiplied, these quantities are often called adjoints
\cite{Griewank1991,Griewank1989,Bischof1992}.
In summary, the result of forward mode AD will be a matrix-vector
product, $\D{\vec{F}}{\vec{u}} \cdot \vec{e}$, for some vector $\vec{e}$,
and the result of reverse mode AD will be a vector-matrix product,
$\vec{e}^T \cdot \D{\vec{F}}{\vec{u}}$.

%In cases where the function output is vector-valued,
%the forward mode can be used to calculate the derivative of the full vector.
%For a vector-valued input, the reverse mode yields the full function gradient.

In order to be able to carry out the derivative calculations in reverse order, reverse
mode requires all operations to be recorded, e.g., in a graph representing the
program structure.
This imposes greater requirements on system memory \cite{Barthelemy1995}.

For the calculation of quadratic Jacobian matrices, the same number of runs of
the derivative function are required, independently of the selected mode.
However, the forward mode promises both an easier implementation and lower
memory requirements under these conditions.

In addition to the different modes for ordering the derivative evaluations,
there are also different methods how to implement them.
The most popular methods are source code transformation
and operator overloading.

Source code transformation methods aim to work with arbitrary function
code that should ideally not require any preparations.
In order to calculate the derivatives, an additional code
transformation step is added between preprocessing and compiling
of the source code \cite{Bischof1992b}.

In contrast, operator overloading methods require the use of special data types
in the function code.
Overloaded versions of all mathematical operations exist for these data types,
such that it is possible to either evaluate a derivative or
accumulate some representation of a computational graph during function
execution \cite{Corliss1993}.
% accumulate
% a computational graph during function execution \cite{Corliss1993}.

In some contexts, source code transformation methods have been shown to be superior
to operator overloading methods, in terms of runtime efficiency
\cite{Tadjouddine2002,Forth2004}.
It also makes a huge difference which programming language is used.
The operator overloading method of AD requires object-oriented programming,
which can be done quite efficiently in a language such as C++. In Fortran,
on the other hand, such concepts have only recently been added to the language,
and are therefore less flexible and their usage may incur performance losses.
The switch from a function code using intrinsic data types to an object type
suitable for AD with operator overloading would also require great changes in
a Fortran code, which may not be easy to automatize.
%In any case, the switch from a function code using intrinsic data types to an object type
%suitable for AD with operator overloading would require much greater changes in
%a Fortran code than, e.g., in a C++ code.

For further details on the different modes and implementation methods, we refer to
\cite{Bischof1992b} and \cite{Corliss1993}.

\subsection{Code preparation for TAPENADE}
\label{sec:software}
\label{sec:s25incomp}

Our FE software uses features from different standards of the Fortran programming
language. This is quite usual for Fortran programs, since most compilers
are backward compatible.
In contrast to these compilers, TAPENADE makes a clear distinction between
different Fortran standards. It allows the selection of either Fortran 77,
90 or 95, and only supports a subset of the available features of each
of these standards.
In addition, our software uses one popular language extension
that is not part of any official Fortran standard, namely Cray pointers.
TAPENADE does not support this extension.
It can be expected that the difficulties arising from these
circumstances are not unique to our software, but will be encountered in most cases
where a large Fortran software with a long development history is to be
used with TAPENADE.

As a consequence, the source code needs to be adjusted in order to
work with TAPENADE. One option to do this would be manual changes
to make the code work with TAPENADE permanently. The amount of effort
that would be necessary to do this is proportional to the size of
the source code. The adjustments would also have to be repeated for
any further software projects that are to be differentiated.

\begin{figure}[t!]
    \centering
    \scalebox{0.8}{
    \begin{tikzpicture}[node distance=0.75cm]
      % text style
      % styles for the flow chart
      \tikzstyle{process} = [rectangle, minimum height=1cm, text centered, text width=5cm, draw=black, fill=white]
      \tikzstyle{script} = [process, fill=RWTHBlueLight]
      \tikzstyle{place} = [diamond, aspect=4, minimum size=2mm, thick, draw, fill=white]
      % line styles
      \tikzstyle{pre}  = [<-,shorten >=1pt,>=stealth,semithick]
      \tikzstyle{mid}  = [-,semithick]
      \tikzstyle{post} = [->,>=stealth,semithick]
      \node (top3) {Source code};
      \node (cpmovement) [process, below = of top3] {Fortran preprocessor} edge[pre] (top3);
      \node (meshmovement) [script, below = of cpmovement] {Preprocess script} edge[pre] (cpmovement);
      \node (nseq) [process, below = of meshmovement] {TAPENADE} edge[pre] (meshmovement);
      \node (velupdate) [coordinate,below = of nseq] {};
      \node[place] (res) [below = of velupdate] {Recursion?} edge[pre] (nseq);
      \node (nseq2) [script, below = of res] {Postprocess script} edge[pre] (res);
      \node (ii1) [script,right = 2.4cm of velupdate] {Midprocess script};
      \path [post,draw] (res) -|  (ii1);
      \path [post,draw]  (ii1) |- (nseq);
      \node (nseq3) [process, below = of nseq2] {Fortran compiler} edge[pre] (nseq2);
      \node (update) [below = of nseq3] {Executable} edge[pre] (nseq3);
      % draw some bounding boxes around the specific loops
      \begin{scope}[on background layer]
        \node [draw=black,rectangle,fit= (ii1) (top3) (update) (nseq2) (res), inner sep = 5pt] (timeloop) {};
      \end{scope}
    \end{tikzpicture}
    }
    \caption{Framework for code transformations}
    \label{fig:s25workflow}
  \end{figure}
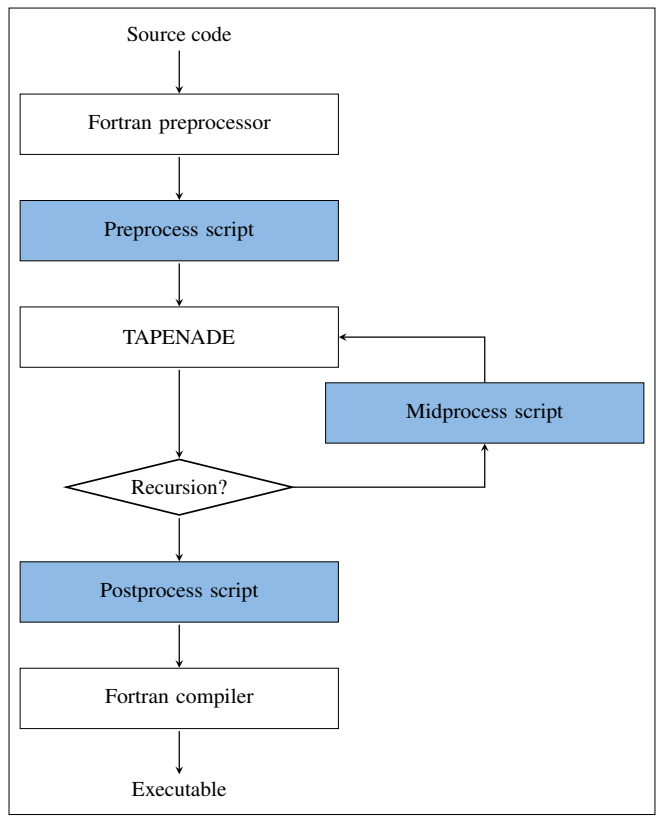

As a more efficient alternative, we initiated the development of a
framework written in Python to automatically prepare a Fortran source code
for use with
TAPENADE. The general structure of the associated workflow
is illustrated in Figure~\ref{fig:s25workflow}. It consists of the
following transformation steps:

\begin{enumerate}
  \item Since TAPENADE needs finished Fortran source code, the first step is
        always to run the Fortran preprocessor. The resulting code is not
        given directly to TAPENADE however, but is first passed to our
        \textit{preprocess script}. Simply put, the job of this script is to
        adjust the source code just enough so that TAPENADE will run
        without errors. Some adjustments can be made without affecting
        the program behavior and can be left in place until compilation.
        In some other cases, the code is changed such that TAPENADE can
        handle it, but that it will no longer be compilable.
  \item In cases where code is differentiated with TAPENADE several times
        recursively,
        some changes are necessary to prepare the TAPENADE output to work
        again as TAPENADE input. This is handled by our \textit{midprocess script}.
  \item As stated before, some changes that are made in the
        \textit{preprocess script} lead to code that will not be accepted
        by the Fortran compiler. These changes need to be reversed by
        the \textit{postprocess script}. Additionally, TAPENADE will
        sometimes produce source code that is also not compilable, so this
        needs to be adjusted in the same script, before the code is passed
        on to the compiler.
\end{enumerate}

In the following, we would like to give a few examples of the transformations
carried out by the scripts. We adjusted the code in such a way that it would
run with the Fortran 95 input mode of TAPENADE version 3.10.

\subsubsection{Cray pointers}

As mentioned before, we use Cray pointers, which are not
standard Fortran and therefore not supported by TAPENADE.
Of course, usage of this construct is not always necessary,
such that one possible option to deal with this issue would
be the complete removal of such pointers from the source code.
However, if the goal is to leave the source code intact,
it is necessary to
preserve Cray pointers during the automatic differentiation process.

In order to do this, we make use of the fact that TAPENADE does
not need correct source code as input. In specific terms, the source code needs to
be syntactically correct, but does not have to be compilable.

The following lines of code provide an example of how a Cray pointer
(called {\tt difptr}) is created as a pointer to an array (called {\tt dif}),
and is subsequently used in the code:

  \begin{lstlisting}
      real* 8 :: dif(nec)
      pointer (difptr, dif)
!     ...
      save difptr
      difptr = malloc(nec*floatsize)
  \end{lstlisting}

The pointer is used like an ordinary variable in the code. Since TAPENADE
does not understand the {\tt pointer} statement, the \textit{preprocess script}
transforms this line to a simple variable declaration:

  \begin{lstlisting}
      real* 8 :: dif(nec)
      real* 8 :: difptr
  \end{lstlisting}

TAPENADE can now understand the source code. If the derivative of the
array {\tt dif} is needed, TAPENADE will add a new array {\tt difd} to
store this:

  \begin{lstlisting}
      real* 8 :: dif(nec)
      real* 8 :: difd(nec)   ! added by TAPENADE
      real* 8 :: difptr
  \end{lstlisting}

The task of our \textit{postprocess script} is then to reverse the
earlier transformation and add another Cray pointer for the new
array {\tt difd}. Any statements involving the original pointer that
relate to memory and allocation will be duplicated for the added
pointer. The resulting source code will take the following form:

  \begin{lstlisting}
      real* 8 :: dif(nec)
      real* 8 :: difd(nec)
      pointer (difptr, dif)
      pointer (difdptr, difd)
!     ...
      save difptr
      save difdptr
      difptr = malloc(nec*floatsize)
      difdptr = malloc(nec*floatsize)
  \end{lstlisting}

\subsubsection{Minimum- and maximum-functions}

Fortran allows minimum- and maximum-functions, such as {\tt dmax1},
to have arbitrary numbers of arguments. TAPENADE, on the other hand, only recognizes
and differentiates these functions if they have exactly two arguments.
The following code excerpt, which is similar to many others in our
software, shows such a function call that TAPENADE cannot handle:

  \begin{lstlisting}
      he(ie) = dmax1(someFunc(a), otherFunc(b, c), d + e, thirdFunc(g) * 2,
     &               h, i)
  \end{lstlisting}

The same functionality can be achieved by using recursive calls to the
same function {\tt dmax1} with just two arguments per function call.
In order to transform the code, our \textit{preprocess script} had to
recognize such functions calls, even if they extend over multiple
lines as in the example. The result of this automatic transformation
looks as follows:

  \begin{lstlisting}
      he(ie) = dmax1(dmax1(dmax1(someFunc(a), otherFunc(b, c)), dmax1(d + e, &
     &thirdFunc(g) * 2)), dmax1(h, i))
  \end{lstlisting}

\subsubsection{Initialization of assumed-size arrays}

Fortran subroutines can have array arguments with assumed size.
When TAPENADE creates derivatives of such arrays, they are defined
the same way, as in the following example:

  \begin{lstlisting}
REAL*8, INTENT(OUT) :: sh2d(lsd, lsd, *)
  \end{lstlisting}

If the changes made to the original array {\tt sh2} depend on some
conditions, TAPENADE tries to automatically initialize the derivative
to zero, using the following statement:

  \begin{lstlisting}
sh2d = 0.0_8
  \end{lstlisting}

Since the size of {\tt sh2d} is unknown at compile time, this statement
leads to an error during compilation. The \textit{postprocess script}
can identify such initializations and remove them from the code. This
is valid for our source code, but different actions may be necessary in
other situations.

\subsubsection{Loop optimization in vector mode}

TAPENADE provides a vector mode that allows the calculation of an arbitrary
number of derivatives simultaneously. This naturally involves the addition of some
loops to the code. TAPENADE creates a new function argument for the number of
derivatives or loop iterations. This argument is usually called {\tt nbdirs}.
Here TAPENADE makes sure to choose a variable name that does not already exist,
as it should, and adds numbers to the name if necessary. This is the desired
behavior in most use cases, but in a special case, this can be improved.

In order to save computational effort, we chose to recursively differentiate
a single function using this vector mode but with respect to different function
arguments. This reduces the number of necessary function calls and could be
shown to improve the performance. As a result, we ended up with a function
containing arguments {\tt nbdirs, nbdirs0, nbdirs1}, etc. The code would then
contain consecutive loops that use some of these variables as the number of
iterations, as shown in the following example:

  \begin{lstlisting}
      do i=1,nbdirs
        tmp2d(i) = tmp1d(i)
      end do
      do i=1,nbdirs0
        tmp3d(i) = tmp1d(i)
      end do
      do i=1,nbdirs1
        tmp4d(i) = tmp1d(i)
      end do
  \end{lstlisting}

All these iteration numbers, {\tt nbdirs, nbdirs0} and {\tt nbdirs1}, are
identical if this is used to calculate Jacobian matrices.
In order to communicate this fact to the compiler,
our \textit{postprocess script} adjusts these loops to use the same variable
{\tt nbdirs} all the time.
This should allow the compiler to recognize that these loops can be merged.
We did, in fact, observe a decrease in computation time following this
change.
%We could observe a decrease in computation time
%following this change, which can probably be attributed to compiler
%optimization.

\subsubsection{Further issues}

There are some other, smaller, issues that have to be handled by the transformation
scripts. Most of these arise from the fact that TAPENADE requires free-form
code when it is run in Fortran 95 mode. They are briefly listed here.

\paragraph{Different comment styles}
Originally, comment lines in Fortran codes had to start with a {\tt c} in
the first column. The newer free-form standard allows comments to begin anywhere but they
have to start with an exclamation mark. Many codes still contain large amounts
of comments in the old style, but TAPENADE only accepts the newer variant when
it is run in Fortran 95 mode. Therefore, the \textit{preprocess script} has to
recognize comments in the old style and transform them.

\paragraph{Multi-line statements}
Similarly to the issue with different comment styles, the free-form standard
has also changed the way how statements that extend over multiple lines are declared.
The original
rule was that an arbitrary character had to appear in the sixth column to mark
a line as a continuation. In free-form code, ampersands before and after continued
code lines are used. TAPENADE once again recognizes only the newer standard, such
that the \textit{preprocess script} needs to transform any multi-line statements
that use the older standard.

\paragraph{Preprocessor lines}
Any lines starting with a hash character are meant for the Fortran preprocessor.
Since macros and include statements definitely have to disappear before TAPENADE is run,
the preprocessor always has to be run beforehand. However, some preprocessors still
leave behind some lines that start with a hash, and TAPENADE complains about this.
The \textit{preprocess script} simply removes all of these lines since they are
not needed.

\paragraph{Statements in newer Fortran standards}
If statements from newer Fortran standards are used, TAPENADE also fails
because it does not recognize them. In some cases, these statements would not
influence the differentiation, so it is possible to remove them temporarily.
For this purpose, our \textit{preprocess script} converts such statements
to comments, and the \textit{postprocess script} reverses this. The only
issue with this process is that TAPENADE sometimes does not preserve all
comments in the code.

\section{Implementation of the Jacobian calculation}
\label{sec:s3impl}

\subsection{Properties of equation systems in the finite element method}
\label{sec:s3implprops}

%\section{Implementation of a PDE in a finite element code}
\label{sec:fe}

In order to discuss the implementation of a solution procedure for a nonlinear partial differential equation
in a finite-element code, we will introduce the following notion of a nonlinear problem:

Consider a function $F:D\subseteq X\to Y'$ that maps a possible subset $D$ of a function space $X$
to the dual space of another function space $Y$, with $X$ and $Y$ both being Banach spaces.
Then the question to seek for a (weak) solution $u\in D$ of a particular nonlinear partial differential equation
can in many cases be cast into the form
\begin{align}
	F(u) &= 0\, ,
\end{align}
or equivalently
\begin{align}
	\label{eqn:genweakform}
	\langle F(u), v\rangle_Y &= 0\quad\forall v\in Y\, .
\end{align}
Here, $\langle \cdot, \cdot\rangle_Y$ denotes the duality bracket, which applies an element of the dual space $Y'$ to
an element of $Y$.

The prototypical example of a PDE, Poisson's equation on a domain $\Omega$,
could in this context be written as $\langle F(u), v\rangle_Y = \int_\Omega \nabla u\cdot\nabla v \, \text{d}\Omega - f(v)$
with a load $f\in Y'$ and the Sobolev spaces $X=Y=H_0^1(\Omega)$ and $Y' = H^{-1}(\Omega)$.

The (conforming) finite element approach of solving such a nonlinear problem consists in choosing
finite dimensional subspaces of $X$ and $Y$, on which a reduced form of the weak problem \eqref{eqn:genweakform}
is solved. For the sake of simplicity, we restrict ourselves in this section
to the standard Galerkin formulation, i.e., we choose equal trial and test spaces $X=Y$,
and denote the finite dimensional subspace by $X_h$.
Selecting a basis $\mathcal{B} = \{e_i\}$ of $X_h$ subsequently gives rise to a nonlinear algebraic equation system
\begin{align}
	\bm{F}(\bm{u}) = \left(\begin{array}{c} \langle F(\sum_i u_i e_i), e_1\rangle_Y \\
	                                        \langle F(\sum_i u_i e_i), e_2\rangle_Y \\ \vdots \end{array}\right)
		&= \bm{0}\quad\mbox{, where } \bm{u} = (u_1,u_2,\ldots)^T \in \mathbb{R}^{\dim(X_h)}\, .
\end{align}

Up to this point we have not made any use of the fact that a PDE, as the name says, is actually composed of
differential operators. One of the key features of differential operators is that they act locally. By that we
mean that considering any function $u\in X$ and some test function $v\in Y$, a perturbation $\delta u$ of $u$
yields
\begin{align}
	\label{eqn:pde_locality}
	\langle F(u+\delta u) - F(u), v\rangle_Y &= 0\, ,
\end{align}
if the essential support of $\delta u$ and $v$ is disjoint. In other words, if we alter the function $u$ in some region
it does not interfere with testing the result $F(u)$ in a separate region.

Considering Newton's method to solve the system of nonlinear algebraic equations, we will furthermore assume that
$\bm{F}$ has a continuous derivative $\D{\vec{F}(\vec{u})}{\vec{u}}$. Given an initial guess $\bm{u}^0$, Newton's
method tries to iteratively approach a solution of the nonlinear equation system by solving the
following sequence of linear equation systems
\begin{align}
    \label{eqn:newton_lineqsys}
      \underbrace{\D{\vec{F}(\vec{u}^k)}{\vec{u}}}_{=: \mat{A}} \, \Delta\vec{u}^{k+1} \; = \; \underbrace{- \vec{F}(\vec{u}^k)}_{=: \vec{b}}\, .
\end{align}

It is especially the construction of this Jacobian of $\bm{F}$, where the locality comes into play. Assuming $e_i$
and $e_j$ are two basis functions with disjoint support, one can readily derive for the particular entries of the
Jacobian matrix
\begin{align}
	A_{ij} &= \lim_{\epsilon\to 0} \frac{1}{\epsilon} \langle F(u+\epsilon\, e_j) - F(u), e_i\rangle_Y = 0\, ,
\end{align}
as follows from \eqref{eqn:pde_locality}. Let $m_{\mathcal{B}}$ be the maximal number of
basis functions that any other basis function in the basis $\mathcal{B}$ can overlap with.
Then, if we select $X_h$ such that $m_{\mathcal{B}} \ll \dim(X_h)$, we can expect $\bm{A}$ to be
quite sparse. The latter of course offers potential savings in computational time and memory if respected during the
construction of the matrix; a savings potential exploited in almost every state-of-the-art code.

More specifically, finite element codes take advantage of this property by splitting the domain $\Omega$ into regions
with a minimal set of non-vanishing basis functions $e_i$. This splitting into elements is usually denoted
as a triangulation $\mathcal{T}$. Furthermore, we let $\eta_e$ denote the ordered index set of non-vanishing
basis functions for a given element $e\in\mathcal{T}$. Considering that we usually deal with Sobolev spaces,
the formation of the matrix $\mat{A}$ and residual vector $\vec{b}$ that make up
this system of equations involves the evaluation of
some integrals on $\Omega$.
The latter allows the finite element code to exploit locality by performing the integral evaluation
element by element, with only a small number $|\eta_e|$ of non-zero basis functions on every element.
More formally, we can split up the function space $Y$ and the function $F$ into pieces that are localized to
one element, i.e., we have for each element
a function space $Y_e$ (with a natural injection $Y\to Y_e$) and a corresponding function $F_e: D_e\to Y_e'$ such that
\begin{align}
	\langle F(u), v \rangle_Y &= \sum_{e\in\mathcal{T}} \langle F_e(u|_e), v|_e\rangle_{Y_e}\, .
\end{align}

Numbering the local basis functions on an element by $\{\tilde{e}_i\}$, we can define
element-level matrices and vectors as
\begin{align}
	\mat{A}_e &:= \D{\vec{F}_e(\vec{u}^k_e)}{\vec{u}_e} \\
	\vec{b}_e &:= - \vec{F}_e(\vec{u}^k_e)\, ,
\end{align}
where $\vec{F}_e$ is given by
\begin{align}
	\bm{F}_e(\bm{u}_e) = \left(\begin{array}{c} \langle F_e(\sum_i u_i \tilde{e}_i), \tilde{e}_1\rangle_{Y_e} \\
	                             \langle F_e(\sum_i u_i \tilde{e}_i), \tilde{e}_2\rangle_{Y_e}
	                             \\ \vdots \end{array}\right)
		\quad\mbox{ with } \bm{u}_e = (u_1,u_2,\ldots)^T \in \mathbb{R}^{|\eta_e|}\, .
\end{align}

  Using the assembly operator, the full matrix and vector $\mat{A}$ and $\vec{b}$
  in \eqref{eqn:newton_lineqsys}
  can then be constructed from the localized element-level variants as

  \begin{equation}
    \begin{aligned}
      \mat{A} \; &= \; \underset{e \in \mathcal{T}}{\boldsymbol{\assembly}} \mat{A}_e
                  := \sum_{e\in\mathcal{T}} \mat{P}_e \mat{A}_e \mat{P}_e^T\, , \\
      \vec{b} \; &= \; \underset{e \in \mathcal{T}}{\boldsymbol{\assembly}} \vec{b}_e
                  := \sum_{e\in\mathcal{T}} \mat{P}_e \vec{b}_e\, .
    \end{aligned}
  \end{equation}
Here, the matrix $\mat{P}_e \in \mathbb{R}^{\dim(X_h) \times |\eta_e|}$ maps back the local element-wise numbering to the global numbering stored in $\eta_e$
and is defined as
\begin{align}
	(\mat{P}_e)_{ij} &:= \begin{cases} 1 & \text{if $(\eta_e)_j = i$,}\\
                                       0 & \text{otherwise.}
                         \end{cases}
\end{align}

  From this perspective, the finite element method
  offers two different points where AD can be applied: It can be used
  to calculate either the element-level matrix $\mat{A}_e$ or the full
  matrix $\mat{A}$. This means that AD is applied either before or
  after the assembly procedure. The two matrices have the following
  properties:

  \begin{itemize}
    \item The \textbf{full matrix} $\mat{A} \in \mathbb{R}^{\dim(X_h) \times \dim(X_h)}$
          grows quadratically
          with the mesh size. For realistic problems, it can get very large
          and it would not be feasible to store it in an uncompressed form.
          Fortunately, as mentioned before, this matrix is also very sparse.
          This means that the matrix can be stored memory-efficiently, e.g.,
          using block
          sparse row (BSR) storage \cite{Smailbegovic2005}.
          In principle, using automatic differentiation for the automatic calculation
          of this matrix would not require any knowledge of the residual calculation
          code, allowing it to be treated as a black box.
          However, methods for sparsity exploitation, which would probably be
          required for this, could profit from some detailed knowledge of the code.

    \item The size of the \textbf{element-level matrix}
          $\mat{A}_e \in \mathbb{R}^{|\eta_e| \times |\eta_e|}$ does not depend
          on the mesh size, but on parameters such as the element type and the
          dimensionality of the differential equation.
          The number of nodes in one element is small for realistic
          cases, such that storage of the matrix is of no concern.
          The structure of this matrix depends on the equations that
          are used, but it is dense in many cases, so it is not
          necessary to exploit matrix sparsity when it is
          accumulated using AD.
  \end{itemize}

  Finite element assembly codes always need to deal with the
  sparsity of the full matrix where memory and computational
  efficiency are concerned. This infrastructure can be reused
  if AD is applied on the element level, removing the need for
  any additional sparsity handling techniques.
  We decided to use this option, since it seems like a much more
  natural approach that is closer to the general idea of finite
  element codes.

\subsection{Efficient calculation of the element-level matrix}
\label{sec:s3jacfe}

  In order to discuss the automatic creation of the code for
  the element-level matrix in more detail, we consider the following
  program function:

  \begin{lstlisting}
      subroutine elemRes(u, res)

      real(8), intent(in)  :: u(ndf*nen)
      real(8), intent(out) :: res(ndf*nen)
!     ...
  \end{lstlisting}

  This function takes an array \texttt{u} for the vector $\bm{u}_e$
  as input, calculates the residual vector $\bm{F}_e\left(\bm{u}_e\right)$,
  and stores it in \texttt{res}.
  We use the variables \texttt{ndf} for the number of degrees of freedom
  and \texttt{nen} for the number of nodes per element, such that
  \texttt{ndf*nen} is the value of $|\eta_e|$. For this to apply, we assume
  equal-order interpolation for all degrees of freedom.
  In practice, a function for the residual calculation would probably
  need additional arguments, but these are not relevant here.

  We now use TAPENADE's vector mode to differentiate this function.
  This yields a new function with a signature similar to the following:

  \begin{lstlisting}
      subroutine elemRes_dv(u, ud, res, resd, nbdirs)

      real(8), intent(in)  :: u(ndf*nen)
      real(8), intent(in)  :: ud(nbdirsmax, ndf*nen)
      real(8), intent(out) :: res(ndf*nen)
      real(8), intent(out) :: resd(nbdirsmax, ndf*nen)
!     ...
  \end{lstlisting}

  The new arguments \texttt{ud} and \texttt{resd} are two-dimensional
  arrays. If we identify the array \texttt{ud} with a
  matrix $\bm{C}$, and \texttt{resd} with
  another matrix $\bm{D}$, where $\bm{C}, \bm{D} \in \mathbb{R}^{|\eta_e|\times|\eta_e|}$,
  the differentiated function will now calculate

  \begin{equation}
    %\bm{D} \; := \; \bm{C} \cdot \left( \D{\vec{F}_e(\vec{u}_e)}{\vec{u}}\right)^T
    %\; = \; \bm{C} \cdot \bm{A}_e^T \\
    \bm{D}^T \; := \; \D{\vec{F}_e(\vec{u}_e)}{\vec{u}_e} \cdot \bm{C}^T
    \; = \; \bm{A}_e \cdot \bm{C}^T
  \end{equation}

  in addition to the residual \texttt{res}\footnote{We also set both
  \texttt{nbdirs} and \texttt{nbdirsmax} to the value of $|\eta_e|$.}.
  By selecting $\bm{C} = \bm{I}$,
  the matrix $\bm{D}$ will take the value of the transpose of the matrix
  $\bm{A}_e$ that we are looking for.

  With the decision to use AD on the element level, we could avoid
  dealing with sparsity in the Jacobian matrix.
  However, in the program function
  \texttt{elemRes\_dv} we now face another sparsity issue. With $\bm{C}$
  being a diagonal matrix, most entries in the array \texttt{ud} will be
  zero. Since the function \texttt{elemRes\_dv} is designed to work
  with arbitrary values of \texttt{ud}, it contains operations to handle
  all array entries. As a result, many operations in the function will
  be multiplications by zero if it is used for Jacobian calculation. This
  in turn renders many intermediate variables in this function redundant,
  as well as many of the operations that they are involved in.

  This is due to the fact that TAPENADE only operates on the level
  of complete arrays and does not allow for distinctions in the
  treatment of individual array entries\footnote{The latter would be
  very difficult to implement whenever variables or expressions are used as array
  indices, which is the rule rather than the exception.}.
  An obvious solution to this problem would
  be the partition of the array \texttt{u} into individual scalar variables
  for all non-zero entries.
  This way, individual functions could be created using AD that
  are optimized for the calculation of one specific column of $\bm{A}_e$.
  Each of these functions would require just one input sensitivity
  for the current scalar variable, instead of a full array with
  $|\eta_e|$ entries. This means
  that the number of input variables for the differentiated function
  would be reduced from $|\eta_e|^2$
  to $|\eta_e|$, with a corresponding reduction in the operation count.
  The only requirement for this procedure would be the knowledge
  of the value of $|\eta_e|$ at compile time, and therefore the necessity to keep
  it fixed.
  Unfortunately, this size is not known, since it depends on the
  element type in use, and the code is designed to be flexible
  in this regard.

  However, even a partition into several smaller arrays, instead
  of scalars, can be expected to improve the situation.
  We may not know the value of $|\eta_e|$ which is represented by
  \texttt{ndf*nen} in the function code, but we do know the number
  of degrees of freedom \texttt{ndf}.
  Assuming this number to be $6$ for now,
  this means that we can split the array \texttt{u} into six smaller
  arrays \texttt{u1} through \texttt{u6} of size \texttt{nen}.
  For each single one of these, it is now possible to create optimized
  differentiated functions, such as:

  \begin{lstlisting}
      subroutine elemRes_dv1(u1, u1d, u2, u3, u4, u5, u6, res, resd, nbdirs)

      real(8), intent(in)  :: u1(nen)
      real(8), intent(in)  :: u1d(nbdirsmax, nen)
      real(8), intent(in)  :: u2(nen)
      real(8), intent(in)  :: u3(nen)
      real(8), intent(in)  :: u4(nen)
      real(8), intent(in)  :: u5(nen)
      real(8), intent(in)  :: u6(nen)
      real(8), intent(out) :: res(ndf*nen)
      real(8), intent(out) :: resd(nbdirsmax, ndf*nen)
!     ...
  \end{lstlisting}

  \texttt{nbdirs} and \texttt{nbdirsmax} now take the lower value of
  \texttt{nen} instead of \texttt{ndf*nen}. The matrix stored in array
  \texttt{u1d} will of course still be a diagonal matrix, but the percentage
  of non-zero entries will be higher than it was in \texttt{ud}.
  In total, the number of input sensitivies would be reduced from
  \texttt{ndf*nen*ndf*nen} entries in \texttt{ud} to \texttt{ndf*nen*nen}
  entries in \texttt{u1d} through \texttt{u6d}. This is a reduction by
  a factor of six.

\begin{figure}[h!]
    \centering
    \scalebox{1.}{
    \begin{tikzpicture}[node distance=0cm]
      \def \si {0.18};
      \foreach \cx in {0, ..., 5}
        \foreach \cy in {0, ..., 5}
          \foreach \x in {0, ..., 2}
            \foreach \y in {0, ..., 2}
              \draw [fill=white, draw=black] (\cx*3*\si+\x*\si,-\cy*3*\si-\y*\si) rectangle (\cx*3*\si+\x*\si+\si,-\cy*3*\si-\y*\si-\si);
      \foreach \cx in {0, ..., 5}
        \foreach \x in {0, ..., 2}
          \draw [fill=black, draw=black] (\cx*3*\si+\x*\si,-\cx*3*\si-\x*\si) rectangle (\cx*3*\si+\x*\si+\si,-\cx*3*\si-\x*\si-\si);
      \foreach \cx in {0, ..., 2}
        \draw [fill=RWTHRed, draw=black] (\cx*6*\si+3*\si,-\cx*6*\si-3*\si) rectangle (\cx*6*\si+3*\si+\si,-\cx*6*\si-3*\si-\si);
      \def \ox {\si*32};
      \def \oy {-\si*0};
      \def \mx {\si*6};
      \def \my {\si*8};
      \foreach \cx in {0, 1}
        \foreach \cy in {0, ..., 2}
          \foreach \x in {0, ..., 2}
            \foreach \y in {0, ..., 2}
              \draw [fill=white, draw=black] (\ox+\cx*\mx+\x*\si,-\oy-\cy*\my-\y*\si) rectangle (\ox+\cx*\mx+\x*\si+\si,-\oy-\cy*\my-\y*\si-\si);
      \foreach \cx in {0, 1}
        \foreach \cy in {0, ..., 2}
          \foreach \x in {0, ..., 2}
            \draw [fill=black, draw=black] (\ox+\cx*\mx+\x*\si,-\oy-\cy*\my-\x*\si) rectangle (\ox+\cx*\mx+\x*\si+\si,-\oy-\cy*\my-\x*\si-\si);
      \foreach \x in {0, ..., 2}
        \draw [fill=RWTHRed, draw=black] (\ox+1*\mx+\x*\si,-\oy-0*\my-\x*\si) rectangle (\ox+1*\mx+\x*\si+\si,-\oy-0*\my-\x*\si-\si);
      \node at (\ox+0*\mx+1.5*\si,-\oy-0*\my+2*\si) {\texttt{u1d}};
      \node at (\ox+0*\mx+1.5*\si,-\oy-1*\my+2*\si) {\texttt{u2d}};
      \node at (\ox+0*\mx+1.5*\si,-\oy-2*\my+2*\si) {\texttt{u3d}};
      \node at (\ox+1*\mx+1.5*\si,-\oy-0*\my+2*\si) {\texttt{u4d}};
      \node at (\ox+1*\mx+1.5*\si,-\oy-1*\my+2*\si) {\texttt{u5d}};
      \node at (\ox+1*\mx+1.5*\si,-\oy-2*\my+2*\si) {\texttt{u6d}};
      \node [scale=3.] at (25*\si,-9*\si) {\MVRightarrow};
      \node at (9*\si,2*\si) {\texttt{ud}};
    \end{tikzpicture}
    }
    \caption{Graphical representation of the sensitivity arrays used as input arguments
             for the different functions created with AD. Filled squares represent non-zero
             values, unfilled squares represent zeros.}
    \label{fig:s3inputdof}
  \end{figure}
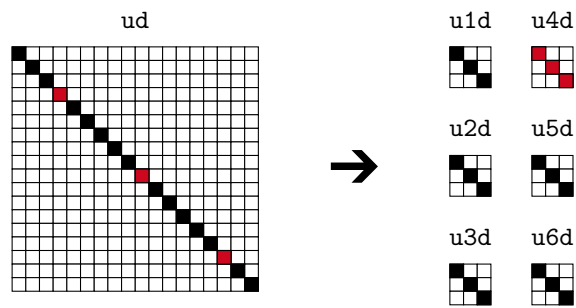

  In Figure~\ref{fig:s3inputdof}, the change is visualized for a set of
  equations with six degrees of freedom and three nodes per element.

\section{Numerical example}
\label{sec:example}

  Now that an approach for the automatic creation of the Jacobian code
  has been described, we will continue by discussing the application of this
  approach on a set of differential equations.
  We will first describe the governing equations and then outline some
  numerical issues that were encountered during the implementation.
  Finally, we will present the results of a test simulation that was run
  using both the automatically created Jacobian code and a manually
  derived variant.

%\subsection{Viscoelastic flow}

  The first step is the choice of the differential equations.
  With the objective of simulating the flow phenomena associated with,
  e.g., blood pumps or plastics injection processes, one new
  set of equations that we are interested in are the viscoelastic fluid equations.
  In particular, we choose the Oldroyd-B constitutive model for this illustration.
  The automatic creation of the Jacobian matrix for these equations can
  be expected to be more challenging than, e.g., for the Navier--Stokes
  equations for Newtonian fluids, since more complicated mathematical
  terms are involved. For this reason, they should be particularly
  suitable for a trial run of our method.

\subsection{Governing equations}

  Viscoelastic flow is also governed by the incompressible Navier--Stokes equations,
  which already involve velocity $\bm{u}\in\mathbb{R}^{n_{sd}}$ and pressure $p$ as
  principal degrees of freedom.
  We extend these with the Oldroyd-B constitutive model in the log-conformation
  formulation \cite{Knechtges2014,Knechtges2015} to account for viscoelastic phenomena.
  This formulation was originally proposed in \cite{Fattal2004}.
  For simplicity, we only consider a two-dimensional simulation domain ($n_{sd}=2$).

  More specifically, we introduce a symmetric tensorial field $\bm{\Psi}$, the so-called log-conformation
  field, which assigns to each point in the fluid domain a symmetric matrix in $\mathbb{R}^{n_{sd}\times n_{sd}}$.
  This field contributes through a matrix exponential function to the polymeric stress tensor
  $\bm{T}_P = \frac{\mu_P}{\lambda} \left(e^{\bm{\Psi}}-\bm{1}\right)$, which then enters the Navier--Stokes
  equations
  \begin{gather}
  \label{eqn:navierstokes_cont}
	\nabla\cdot \bm{u} = 0\\
  \label{eqn:navierstokes_mom}
	\rho (\partial_t + \bm{u}\cdot\nabla) \bm{u} + \nabla p
	    - \mu_S \nabla\cdot\left(\nabla\bm{u}+\nabla\bm{u}^T\right)
	    - \nabla \cdot \bm{T}_P = \bm{0}\, .
  \end{gather}
  Here, $\mu_P$ and $\mu_S$ denote the polymeric and solvent viscosity, respectively.
  $\lambda$ is the relaxation time and $\rho$ the density. Moreover, we need a
  constitutive equation to close the formulation. In our case, this equation is for the specific choice
  of the Oldroyd-B model given by
  \begin{align}
  \begin{split}
		&\partial_t\bm{\Psi} + (\bm{u}\cdot \nabla) \bm{\Psi} + [\bm{\Psi},\Omega(\bm{u})]
				+ \frac{1}{\lambda} \left(\bm{1}-e^{-\bm{\Psi}}\right) - 2\varepsilon(\bm{u})\\ &\quad
				- 2 \left(\begin{array}{cc} -\bm{\Psi}_{12} & \gamma(\bm{\Psi}) \\
				                     \gamma(\bm{\Psi}) & \bm{\Psi}_{12}\end{array}\right)
				\left[\gamma(\bm{\Psi})\varepsilon(\bm{u})_{12}
				     - \bm{\Psi}_{12}\gamma(\varepsilon(\bm{u}))\right] \cdot f(\bm{\Psi}) =
				\bm{0}\, ,
  \end{split}
  \label{eqn:2DPsi}
  \end{align}
  where $\varepsilon(\bm{u}) = \frac{1}{2}(\nabla\bm{u}+\nabla\bm{u}^T)$ denotes the rate of strain tensor
  and $\Omega(\bm{u}) = \frac{1}{2}(\nabla\bm{u}-\nabla\bm{u}^T)$ the vorticity tensor. Furthermore,
  $\gamma(\bm{M})$ is an abbreviation for $\frac{1}{2}(\bm{M}_{11}-\bm{M}_{22})$, whereas $f$ is defined by
  \begin{align*}
      f(\bm{\Psi}) =& \frac{1}{\gamma(\bm{\Psi})^2+\bm{\Psi}_{12}^2} \left(
			\sqrt{\gamma(\bm{\Psi})^2+\bm{\Psi}_{12}^2} + \frac{2\sqrt{\gamma(\bm{\Psi})^2+\bm{\Psi}_{12}^2}}{\exp\left(2\sqrt{\gamma(\bm{\Psi})^2+\bm{\Psi}_{12}^2}\right)-1}
			- 1 \right)\, .
  \end{align*}
  For the two-dimensional case, this leads to a total of six
  degrees of
  freedom---$\bm{\Psi}_{11},\bm{\Psi}_{12},\bm{\Psi}_{22},\bm{u}_1,\bm{u}_2,p$---,
  which are governed by the
  Eqs.~\eqref{eqn:navierstokes_cont}-\eqref{eqn:2DPsi}.

  We will not go much further into
  the details of discretization, for which we refer to
  \cite{Knechtges2014}.
  It is however important to note that we obtain a non-linear algebraic equation system
  which is suitable for the Newton--Raphson algorithm.

  Looking at the set of equations, and regardless of the discretization that is used, it becomes apparent that
  one also needs to choose a numerical approximation of the matrix exponential function.
  Following the implementation in \cite{Knechtges2014}, the results in this paper use a
  scaling/squaring approach \cite{Moler2003}. The main principle behind this approach is to use
  the identity $\exp(\bm{M}) = \left(\exp(\bm{M}/2^j)\right)^{2^j}$, where $j$ is chosen such that $\exp(\bm{M}/2^j)$ is
  sufficiently well approximated by a Padé-approximant $R_{6,6}(\bm{M}/2^j)$.
  More precisely, in pseudocode this may look like the following algorithm:
  \begin{spacing}{1.1}
  \begin{algorithmic}[1]
  \Procedure{expm}{$\bm{M}$}
  \State $infnorm\gets ||\bm{M}||_{\infty}$
  \If{$infnorm\geq 10^{-12}$}
    \State $j\gets \max(1,\lceil\log_2(infnorm)\rceil)$
  \Else
    \State \textbf{return} $\bm{1}$\label{alg:prelimreturn}
  \EndIf
  \State $\bm{M}\gets 2^{-j}\bm{M}$
  \State $\bm{L}\gets \bm{1} - \frac{1}{2}\bm{M} + \frac{5}{44}\bm{M}^2 - \frac{1}{66}\bm{M}^3
                             + \frac{1}{792}\bm{M}^4 - \frac{1}{15840}\bm{M}^5 + \frac{1}{665280}\bm{M}^6$
  \State $\bm{K}\gets \bm{1} + \frac{1}{2}\bm{M} + \frac{5}{44}\bm{M}^2 + \frac{1}{66}\bm{M}^3
                             + \frac{1}{792}\bm{M}^4 + \frac{1}{15840}\bm{M}^5 + \frac{1}{665280}\bm{M}^6$
  \State $\bm{K}\gets \bm{L}^{-1}\bm{K}$
  \State $i\gets 1$
  \While{$i \leq j$}
    \State $\bm{K}\gets\bm{K}^2$
    \State $i\gets i+1$
  \EndWhile
  \State \textbf{return} $\bm{K}$
  \EndProcedure
  \end{algorithmic}
  \end{spacing}

\subsection{Numerical issues}

  When it now comes to applying AD to the \textsc{expm} algorithm, one issue
  becomes directly apparent:
  When encountering conditional statements, AD (as implemented by TAPENADE)
  will always differentiate each
  branch independently. This means that the derivative that is calculated will
  only depend on the single branch that is executed for the current evaluation
  point. It is up to the user to ensure that the derivatives of the different
  branches of an if-block match.
  In mathematical terms, if-conditions can also cause discontinuities in the
  function. If this happens, jumps in the derivative of the function can occur.
  The user must try to avoid this or otherwise make sure that these jumps
  remain so small that they can be neglected.
  In the context of our viscoelastic fluid flow simulation this turned out to be
  a crucial point.
  %When it now comes to applying AD to the \textsc{expm} algorithm, one drawback of AD becomes directly apparent:
  %If-conditions, which in mathematical terms could cause a discontinuity in the function subject to differentiation,
  %are simply ignored by the AD algorithm. In the context of our viscoelastic fluid flow simulation this turned out to be
  %quite crucial.

  More specifically, the return statement in line \ref{alg:prelimreturn} leads to a zero derivative of $e^{\bm{M}}$
  as long as $\bm{M}$ is sufficiently close to zero. Considering now Eq.~\eqref{eqn:2DPsi} in the stationary limit, i.e. $\partial_t \bm{\Psi}=0$,
  and starting from an initial guess $\bm{u}=\bm{0},\bm{\Psi}=\bm{0}$, the only term in Eq.~\eqref{eqn:2DPsi} that still contributes to the left-hand side of the
  linear equation system is the derivative of $\frac{1}{\lambda} e^{-\bm{\Psi}}$ with respect to $\bm{\Psi}$. But since the latter
  evaluates to zero, these rows of the final matrix will be zero, leading to an unsolvable equation system.

  A remedy to this issue is of course quite easy: One needs to augment the return statement in line \ref{alg:prelimreturn} by the next order term in
  the Taylor series. In other words returning $\bm{1}+\bm{M}$ solves the problem.
  Alternatively, one may just set $j=0$ without returning from the function.

  This example quickly makes clear that AD is not always a panacea, which can be applied in a black-box fashion to the software.
  There are situations where the discrete nature of computations shines through, and tracing
  and fixing these issues still requires a deep understanding of the underlying code. This becomes even more pronounced, when one starts looking at other
  implementations of the matrix exponential function, which shall be mentioned here for completeness.

  Since we are dealing with a symmetric matrix $\bm{\Psi}$, one could also consider the implementation of the matrix exponential function using
  the spectral decomposition of $\bm{\Psi}$ into its eigenvalues and eigenvectors. We introduce $\lambda_i$ and $\bm{e}_i$, which are
  for $1\leq i\leq n_{sd}$ the eigenvalues and eigenvectors, respectively, for a general symmetric matrix $\bm{M}$. The matrix exponential
  function is then given by
  \begin{align}
      e^{\bm{M}} =& \sum_{i=1}^{n_{sd}} e^{\lambda_i} \bm{P}_i\, ,
  \end{align}
  where $\bm{P}_i = \bm{e}_i\bm{e}_i^T$ is the projection operator onto the eigenspace corresponding to $\lambda_i$.
  Using this definition as a basis for the implementation, a typical AD algorithm would also yield suboptimal results.
  This can in particular be attributed to the fact that for degenerate eigenvalues, i.e.,
  eigenvalues with multiplicity greater than one, the choice of eigenvectors that span the corresponding eigenspace is ambiguous. In this light, a small perturbation of the
  matrix that draws two equal eigenvalues apart may also yield an instantaneous jump-like reaction of the corresponding projections $\bm{P}_i$. The latter means that
  a total derivative of $\bm{P}_i$ at this point does not need to exist\footnote{For more details on this issue in the view of perturbation theory
  \cite[\S II.5.8]{Kato}.}. It should be stressed that this is not an artifact of finite precision arithmetic, but inherent
  to the problem formulation. As such, an AD-derivative of $\bm{P}_i$, and even the eigenvalues $\lambda_i$,
  would potentially suffer from the same issue.
  In the end, all the singularities in the derivatives of the projection operators
  have to cancel, since the derivative of the matrix exponential function is well-defined. In two dimensions, the latter can be seen from the following closed formulation
  \cite{Knechtges2014}
  \begin{align}
	\partial_i e^{\bm{M}(x)} &= e^{\bm{M}/2}\left(\partial_i \bm{M}(x) + \left(\begin{array}{cc} - \bm{M}_{12} & \gamma(\bm{M})\\
				\gamma(\bm{M}) & \bm{M}_{12}\end{array}\right)\left[\gamma(\bm{M})\partial_i \bm{M}_{12} - \bm{M}_{12} \partial_i \gamma(\bm{M})\right]
				\cdot g(\bm{M})\right)e^{\bm{M}/2}\, ,
	\label{eqn:expm_manuel}
  \end{align}
  where
  \begin{align*}
	g(\bm{M}) =& \left(\gamma(\bm{M})^2+\bm{M}_{12}^2\right)^{-3/2}\cdot\left(\sinh\left(\sqrt{\gamma(\bm{M})^2+\bm{M}_{12}^2}\right) -
				\sqrt{\gamma(\bm{M})^2+\bm{M}_{12}^2}\right)\, .
  \end{align*}
  A typical AD algorithm is of course not aware of these cancellations that lead to a removable singularity. It should not be left unmentioned that there is ongoing research
  to automatically determine removable singularities by means of one-sided Laurent expansions \cite{Griewank2008}.

  The principal conclusion here is that one always needs to be aware that an AD
algorithm basically has no other choice but to differentiate the exact
operations that it finds in the program code.
The derivative of a numerical code that just approximates the result of some
mathematical function may actually not be a good approximation of the true
function derivative. This is independent of how accurate the result of the
numerical code may be.
In practice, this means that when AD is applied to an unknown numerical code
that is known to produce results that are accurate enough for some particular purpose,
there is no guarantee that the resulting derivative is of any use.

\subsection{Simulation setup}

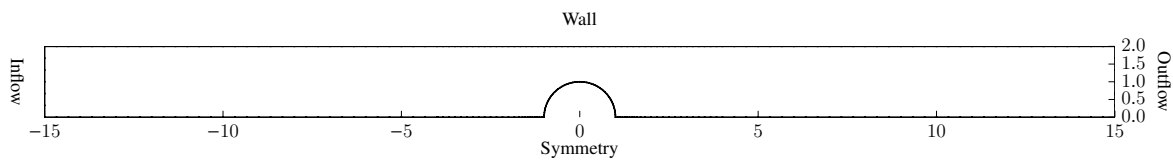
\begin{figure}[h!]
    \centering
    \scalebox{1.}{
    \begin{tikzpicture}[node distance=0.75cm]
      % text style
      % styles for the flow chart
      \tikzstyle{desc} = [text centered,scale=.7]
      \node (pic) {\includegraphics[width=15cm]{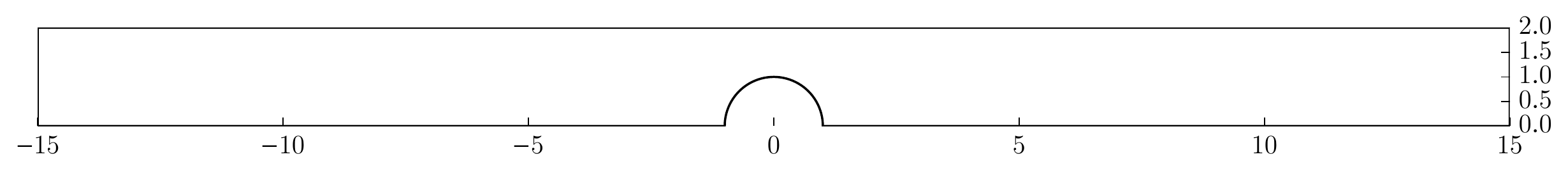}};
      \coordinate (nw) at ($(pic.north west) + (0.5cm,-0.4cm)$);
      \coordinate (sw) at ($(pic.south west) + (0.5cm,0.6cm)$);
      \coordinate (ne) at ($(pic.north east) + (-0.7cm,-0.4cm)$);
      \coordinate (se) at ($(pic.south east) + (-0.7cm,0.6cm)$);
      \coordinate (nc) at ($0.5*(nw)+0.5*(ne)$);
      \coordinate (sc) at ($0.5*(sw)+0.5*(se)$);
      \coordinate (wc) at ($0.5*(nw)+0.5*(sw)$);
      \coordinate (ec) at ($0.5*(ne)+0.5*(se)$);
      \node (sym) [desc,below = 0.25cm of sc] {Symmetry};
      \node (wall) [desc,above = 0.2cm of nc] {Wall};
      \node (in) [desc,left = 0.2cm of wc, rotate=-90,anchor=north] {Inflow};
      \node (out) [desc,right = 0.5cm of ec,rotate=-90,anchor=south] {Outflow};
    \end{tikzpicture}
    }
    \caption{Test case setup for flow around a cylinder.}
    \label{fig:modelcylsetup}
  \end{figure}

  As a test case for the comparison of the automatically created
  Jacobian code with the manually derived variant,
  we performed the simulation of a stationary flow around a cylinder.
  The test case setup is shown in Figure~\ref{fig:modelcylsetup}.

  The simulation domain is cut in half to simplify the setup,
  with symmetry boundary conditions prescribed on the lower
  boundary.
  The flow inlet is located on the left-hand side boundary,
  the outlet on the right-hand side. The top boundary and
  cylinder are treated as walls, using no-slip boundary
  conditions. For the inlet conditions, we use a parabolic
  velocity profile, with a maximum velocity value of $3$.
  We use constant values of density ($\rho = 1$), polymer viscosity
  ($\mu_p = 0.41$), and solvent viscosity ($\mu_s = 0.59$). Only the
  relaxation time $\lambda$ is varied, with values chosen between $0.05$
  and $0.45$ to achieve different Weissenberg numbers.
  The interpolation spaces for all solution fields are
  of second order ($\mathbb{P}_2$).
  The test case was originally used in \cite{Knechtges2014},
  where it is explained in greater detail.

\subsection{Results}
\label{sec:s35res}

  Convergence plots of simulations with two different mesh resolutions are shown
  in Figures~\ref{fig:s35newtonits48}~and~\ref{fig:s35newtonits96}.
  The quadratic convergence behavior that was attained using the manual derivation
  of the Jacobian code (called
  ``no AD'' in the plot) can
  also be reproduced with the matrix created by automatic differentiation (``AD'').

  %\begin{figure}[h!]
  %  \centering
  %  \begin{minipage}{0.48\textwidth}
  %    \input{../Post/paper/conv/res/newtonits48.tex}
  %    \caption{Convergence of the residual with different \mbox{implementations} and mesh ``$48$''.}
  %    \label{fig:s35newtonits48old}
  %  \end{minipage}
  %  \hfill
  %  \begin{minipage}{0.48\textwidth}
  %    \input{../Post/paper/conv/res/newtonits96.tex}
  %    \caption{Convergence of the residual with different \mbox{implementations} and mesh ``$96$''.}
  %    \label{fig:s35newtonits96old}
  %  \end{minipage}
  %\end{figure}

  \begin{figure}[h!]
    \centering
    \begin{minipage}{0.48\textwidth}
      \begin{tikzpicture}[gnuplot]
%% generated with GNUPLOT 5.0p0 (Lua 5.2; terminal rev. 99, script rev. 100)
%% Tue Jul 19 11:18:25 2016
\tikzset{every node/.append style={scale=0.75}}
\path (0.000,0.000) rectangle (8.000,7.000);
\gpcolor{color=gp lt color border}
\gpsetlinetype{gp lt border}
\gpsetdashtype{gp dt solid}
\gpsetlinewidth{1.00}
\draw[gp path] (1.266,0.739)--(1.446,0.739);
\draw[gp path] (7.585,0.739)--(7.405,0.739);
\node[gp node right] at (1.128,0.739) {1e-16};
\draw[gp path] (1.266,1.166)--(1.356,1.166);
\draw[gp path] (7.585,1.166)--(7.495,1.166);
\draw[gp path] (1.266,1.594)--(1.446,1.594);
\draw[gp path] (7.585,1.594)--(7.405,1.594);
\node[gp node right] at (1.128,1.594) {1e-14};
\draw[gp path] (1.266,2.021)--(1.356,2.021);
\draw[gp path] (7.585,2.021)--(7.495,2.021);
\draw[gp path] (1.266,2.449)--(1.446,2.449);
\draw[gp path] (7.585,2.449)--(7.405,2.449);
\node[gp node right] at (1.128,2.449) {1e-12};
\draw[gp path] (1.266,2.876)--(1.356,2.876);
\draw[gp path] (7.585,2.876)--(7.495,2.876);
\draw[gp path] (1.266,3.304)--(1.446,3.304);
\draw[gp path] (7.585,3.304)--(7.405,3.304);
\node[gp node right] at (1.128,3.304) {1e-10};
\draw[gp path] (1.266,3.731)--(1.356,3.731);
\draw[gp path] (7.585,3.731)--(7.495,3.731);
\draw[gp path] (1.266,4.158)--(1.446,4.158);
\draw[gp path] (7.585,4.158)--(7.405,4.158);
\node[gp node right] at (1.128,4.158) {1e-08};
\draw[gp path] (1.266,4.586)--(1.356,4.586);
\draw[gp path] (7.585,4.586)--(7.495,4.586);
\draw[gp path] (1.266,5.013)--(1.446,5.013);
\draw[gp path] (7.585,5.013)--(7.405,5.013);
\node[gp node right] at (1.128,5.013) {1e-06};
\draw[gp path] (1.266,5.441)--(1.356,5.441);
\draw[gp path] (7.585,5.441)--(7.495,5.441);
\draw[gp path] (1.266,5.868)--(1.446,5.868);
\draw[gp path] (7.585,5.868)--(7.405,5.868);
\node[gp node right] at (1.128,5.868) {1e-04};
\draw[gp path] (1.266,6.296)--(1.356,6.296);
\draw[gp path] (7.585,6.296)--(7.495,6.296);
\draw[gp path] (1.266,6.723)--(1.446,6.723);
\draw[gp path] (7.585,6.723)--(7.405,6.723);
\node[gp node right] at (1.128,6.723) {1e-02};
\draw[gp path] (1.266,0.739)--(1.266,0.919);
\draw[gp path] (1.266,6.723)--(1.266,6.543);
\node[gp node center] at (1.266,0.508) {$0$};
\draw[gp path] (2.319,0.739)--(2.319,0.919);
\draw[gp path] (2.319,6.723)--(2.319,6.543);
\node[gp node center] at (2.319,0.508) {$2$};
\draw[gp path] (3.372,0.739)--(3.372,0.919);
\draw[gp path] (3.372,6.723)--(3.372,6.543);
\node[gp node center] at (3.372,0.508) {$4$};
\draw[gp path] (4.426,0.739)--(4.426,0.919);
\draw[gp path] (4.426,6.723)--(4.426,6.543);
\node[gp node center] at (4.426,0.508) {$6$};
\draw[gp path] (5.479,0.739)--(5.479,0.919);
\draw[gp path] (5.479,6.723)--(5.479,6.543);
\node[gp node center] at (5.479,0.508) {$8$};
\draw[gp path] (6.532,0.739)--(6.532,0.919);
\draw[gp path] (6.532,6.723)--(6.532,6.543);
\node[gp node center] at (6.532,0.508) {$10$};
\draw[gp path] (7.585,0.739)--(7.585,0.919);
\draw[gp path] (7.585,6.723)--(7.585,6.543);
\node[gp node center] at (7.585,0.508) {$12$};
\draw[gp path] (1.266,6.723)--(1.266,0.739)--(7.585,0.739)--(7.585,6.723)--cycle;
\node[gp node center,rotate=-270] at (0.184,3.731) {Residual};
\node[gp node center] at (4.425,0.162) {Iteration};
\node[gp node right] at (6.439,6.385) {no AD};
\gpcolor{rgb color={0.800,0.027,0.118}}
\draw[gp path] (6.577,6.385)--(7.309,6.385);
\draw[gp path] (1.793,6.685)--(2.319,6.379)--(2.846,6.402)--(3.372,6.356)--(3.899,6.127)%
  --(4.426,5.696)--(4.952,4.951)--(5.479,3.842)--(6.005,3.146)--(6.532,2.657)--(7.058,2.176);
\gpsetpointsize{4.00}
\gppoint{gp mark 1}{(1.793,6.685)}
\gppoint{gp mark 1}{(2.319,6.379)}
\gppoint{gp mark 1}{(2.846,6.402)}
\gppoint{gp mark 1}{(3.372,6.356)}
\gppoint{gp mark 1}{(3.899,6.127)}
\gppoint{gp mark 1}{(4.426,5.696)}
\gppoint{gp mark 1}{(4.952,4.951)}
\gppoint{gp mark 1}{(5.479,3.842)}
\gppoint{gp mark 1}{(6.005,3.146)}
\gppoint{gp mark 1}{(6.532,2.657)}
\gppoint{gp mark 1}{(7.058,2.176)}
\gppoint{gp mark 1}{(6.943,6.385)}
\gpcolor{color=gp lt color border}
\node[gp node right] at (6.439,6.070) {AD};
\gpcolor{rgb color={0.000,0.329,0.624}}
\draw[gp path] (6.577,6.070)--(7.309,6.070);
\draw[gp path] (1.793,6.685)--(2.319,6.379)--(2.846,6.404)--(3.372,6.359)--(3.899,6.186)%
  --(4.426,5.934)--(4.952,5.421)--(5.479,4.341)--(6.005,2.727)--(6.532,1.832)--(7.058,1.349);
\gppoint{gp mark 2}{(1.793,6.685)}
\gppoint{gp mark 2}{(2.319,6.379)}
\gppoint{gp mark 2}{(2.846,6.404)}
\gppoint{gp mark 2}{(3.372,6.359)}
\gppoint{gp mark 2}{(3.899,6.186)}
\gppoint{gp mark 2}{(4.426,5.934)}
\gppoint{gp mark 2}{(4.952,5.421)}
\gppoint{gp mark 2}{(5.479,4.341)}
\gppoint{gp mark 2}{(6.005,2.727)}
\gppoint{gp mark 2}{(6.532,1.832)}
\gppoint{gp mark 2}{(7.058,1.349)}
\gppoint{gp mark 2}{(6.943,6.070)}
\gpcolor{color=gp lt color border}
\draw[gp path] (1.266,6.723)--(1.266,0.739)--(7.585,0.739)--(7.585,6.723)--cycle;
%% coordinates of the plot area
\gpdefrectangularnode{gp plot 1}{\pgfpoint{1.266cm}{0.739cm}}{\pgfpoint{7.585cm}{6.723cm}}
\end{tikzpicture}
%% gnuplot variables
      \caption{Convergence of the residual with different \mbox{implementations} and mesh ``$48$''.}
      \label{fig:s35newtonits48}
    \end{minipage}
    \hfill
    \begin{minipage}{0.48\textwidth}
      \begin{tikzpicture}[gnuplot]
%% generated with GNUPLOT 5.0p0 (Lua 5.2; terminal rev. 99, script rev. 100)
%% Tue Jul 19 11:18:25 2016
\tikzset{every node/.append style={scale=0.75}}
\path (0.000,0.000) rectangle (8.000,7.000);
\gpcolor{color=gp lt color border}
\gpsetlinetype{gp lt border}
\gpsetdashtype{gp dt solid}
\gpsetlinewidth{1.00}
\draw[gp path] (1.266,0.739)--(1.446,0.739);
\draw[gp path] (7.585,0.739)--(7.405,0.739);
\node[gp node right] at (1.128,0.739) {1e-16};
\draw[gp path] (1.266,1.166)--(1.356,1.166);
\draw[gp path] (7.585,1.166)--(7.495,1.166);
\draw[gp path] (1.266,1.594)--(1.446,1.594);
\draw[gp path] (7.585,1.594)--(7.405,1.594);
\node[gp node right] at (1.128,1.594) {1e-14};
\draw[gp path] (1.266,2.021)--(1.356,2.021);
\draw[gp path] (7.585,2.021)--(7.495,2.021);
\draw[gp path] (1.266,2.449)--(1.446,2.449);
\draw[gp path] (7.585,2.449)--(7.405,2.449);
\node[gp node right] at (1.128,2.449) {1e-12};
\draw[gp path] (1.266,2.876)--(1.356,2.876);
\draw[gp path] (7.585,2.876)--(7.495,2.876);
\draw[gp path] (1.266,3.304)--(1.446,3.304);
\draw[gp path] (7.585,3.304)--(7.405,3.304);
\node[gp node right] at (1.128,3.304) {1e-10};
\draw[gp path] (1.266,3.731)--(1.356,3.731);
\draw[gp path] (7.585,3.731)--(7.495,3.731);
\draw[gp path] (1.266,4.158)--(1.446,4.158);
\draw[gp path] (7.585,4.158)--(7.405,4.158);
\node[gp node right] at (1.128,4.158) {1e-08};
\draw[gp path] (1.266,4.586)--(1.356,4.586);
\draw[gp path] (7.585,4.586)--(7.495,4.586);
\draw[gp path] (1.266,5.013)--(1.446,5.013);
\draw[gp path] (7.585,5.013)--(7.405,5.013);
\node[gp node right] at (1.128,5.013) {1e-06};
\draw[gp path] (1.266,5.441)--(1.356,5.441);
\draw[gp path] (7.585,5.441)--(7.495,5.441);
\draw[gp path] (1.266,5.868)--(1.446,5.868);
\draw[gp path] (7.585,5.868)--(7.405,5.868);
\node[gp node right] at (1.128,5.868) {1e-04};
\draw[gp path] (1.266,6.296)--(1.356,6.296);
\draw[gp path] (7.585,6.296)--(7.495,6.296);
\draw[gp path] (1.266,6.723)--(1.446,6.723);
\draw[gp path] (7.585,6.723)--(7.405,6.723);
\node[gp node right] at (1.128,6.723) {1e-02};
\draw[gp path] (1.266,0.739)--(1.266,0.919);
\draw[gp path] (1.266,6.723)--(1.266,6.543);
\node[gp node center] at (1.266,0.508) {$0$};
\draw[gp path] (2.319,0.739)--(2.319,0.919);
\draw[gp path] (2.319,6.723)--(2.319,6.543);
\node[gp node center] at (2.319,0.508) {$2$};
\draw[gp path] (3.372,0.739)--(3.372,0.919);
\draw[gp path] (3.372,6.723)--(3.372,6.543);
\node[gp node center] at (3.372,0.508) {$4$};
\draw[gp path] (4.426,0.739)--(4.426,0.919);
\draw[gp path] (4.426,6.723)--(4.426,6.543);
\node[gp node center] at (4.426,0.508) {$6$};
\draw[gp path] (5.479,0.739)--(5.479,0.919);
\draw[gp path] (5.479,6.723)--(5.479,6.543);
\node[gp node center] at (5.479,0.508) {$8$};
\draw[gp path] (6.532,0.739)--(6.532,0.919);
\draw[gp path] (6.532,6.723)--(6.532,6.543);
\node[gp node center] at (6.532,0.508) {$10$};
\draw[gp path] (7.585,0.739)--(7.585,0.919);
\draw[gp path] (7.585,6.723)--(7.585,6.543);
\node[gp node center] at (7.585,0.508) {$12$};
\draw[gp path] (1.266,6.723)--(1.266,0.739)--(7.585,0.739)--(7.585,6.723)--cycle;
\node[gp node center,rotate=-270] at (0.184,3.731) {Residual};
\node[gp node center] at (4.425,0.162) {Iteration};
\node[gp node right] at (6.439,6.385) {no AD};
\gpcolor{rgb color={0.800,0.027,0.118}}
\draw[gp path] (6.577,6.385)--(7.309,6.385);
\draw[gp path] (1.793,6.446)--(2.319,6.147)--(2.846,6.169)--(3.372,6.147)--(3.899,5.901)%
  --(4.426,5.343)--(4.952,4.286)--(5.479,2.641)--(6.005,1.752);
\gpsetpointsize{4.00}
\gppoint{gp mark 1}{(1.793,6.446)}
\gppoint{gp mark 1}{(2.319,6.147)}
\gppoint{gp mark 1}{(2.846,6.169)}
\gppoint{gp mark 1}{(3.372,6.147)}
\gppoint{gp mark 1}{(3.899,5.901)}
\gppoint{gp mark 1}{(4.426,5.343)}
\gppoint{gp mark 1}{(4.952,4.286)}
\gppoint{gp mark 1}{(5.479,2.641)}
\gppoint{gp mark 1}{(6.005,1.752)}
\gppoint{gp mark 1}{(6.943,6.385)}
\gpcolor{color=gp lt color border}
\node[gp node right] at (6.439,6.070) {AD};
\gpcolor{rgb color={0.000,0.329,0.624}}
\draw[gp path] (6.577,6.070)--(7.309,6.070);
\draw[gp path] (1.793,6.446)--(2.319,6.147)--(2.846,6.171)--(3.372,6.149)--(3.899,5.903)%
  --(4.426,5.357)--(4.952,4.371)--(5.479,2.619)--(6.005,1.431);
\gppoint{gp mark 2}{(1.793,6.446)}
\gppoint{gp mark 2}{(2.319,6.147)}
\gppoint{gp mark 2}{(2.846,6.171)}
\gppoint{gp mark 2}{(3.372,6.149)}
\gppoint{gp mark 2}{(3.899,5.903)}
\gppoint{gp mark 2}{(4.426,5.357)}
\gppoint{gp mark 2}{(4.952,4.371)}
\gppoint{gp mark 2}{(5.479,2.619)}
\gppoint{gp mark 2}{(6.005,1.431)}
\gppoint{gp mark 2}{(6.943,6.070)}
\gpcolor{color=gp lt color border}
\draw[gp path] (1.266,6.723)--(1.266,0.739)--(7.585,0.739)--(7.585,6.723)--cycle;
%% coordinates of the plot area
\gpdefrectangularnode{gp plot 1}{\pgfpoint{1.266cm}{0.739cm}}{\pgfpoint{7.585cm}{6.723cm}}
\end{tikzpicture}
%% gnuplot variables
      \caption{Convergence of the residual with different \mbox{implementations} and mesh ``$96$''.}
      \label{fig:s35newtonits96}
    \end{minipage}
  \end{figure}

  The two simulations considered here use different meshes. The rougher mesh, called
  ``$48$'', has 2532 elements, the finer mesh ``$96$'' has 10104. We define
  the dimensionless Weissenberg number as $Wi = \lambda \bar{u} / R$, with relaxation time $\lambda$,
  mean inflow velocity $\bar{u}$ and circle radius $R$. We use
  $Wi = 0.8$ in both cases.

  The plots show slight differences between the curves obtained with the
  different implementations. These are due to the fact that the derivatives of
  certain parts of the stabilization term are not considered in the manual
  implementation for reasons concerning the convergence.
  The residual always converges to zero in the end. Assuming a unique
  solution to the differential equation, this implies that the solution fields
  converge to the same values in both cases.

\begin{figure}[h!]
    \centering
    % 9.77 x 1.15, width=15cm, height=1.7656
    \begin{minipage}[b]{15.281cm}
      \includegraphics[width=15cm]{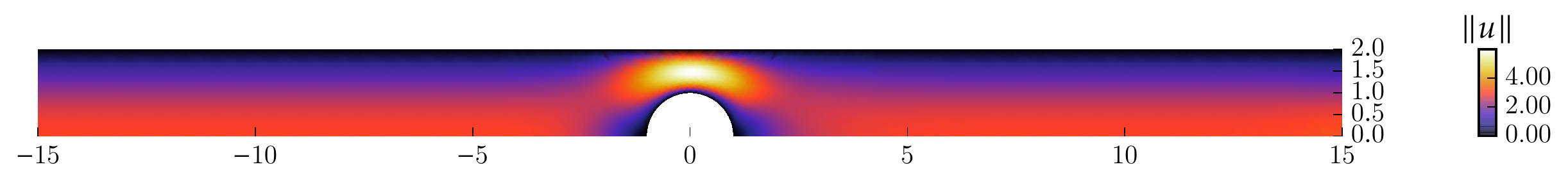}
    \end{minipage}
    \caption{Flow field for mesh ``$48$'' with manually implemented matrix.}
    \label{fig:compadmag48NOAD}
  \end{figure}

\begin{figure}[h!]
    \centering
    % 9.78 x 1.13, width=15.281cm, height=1.7656
    \begin{minipage}[b]{15.281cm}
      \includegraphics[width=15.281cm]{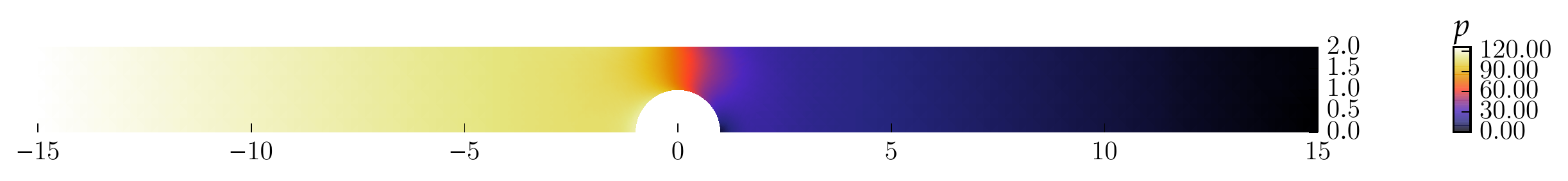}
    \end{minipage}
    \caption{Pressure field for mesh ``$48$'' with manually implemented matrix.}
    \label{fig:compadp48NOAD}
  \end{figure}

\begin{figure}[h!]
    \centering
    % 9.77 x 1.15, width=15cm, height=1.7656
    \begin{minipage}[b]{15.281cm}
      \includegraphics[width=15cm]{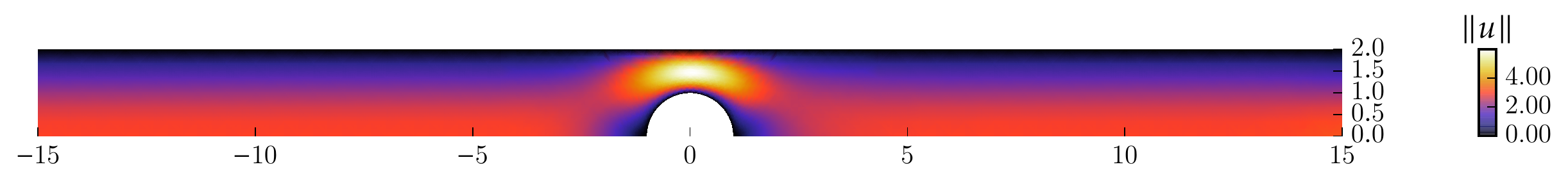}
    \end{minipage}
    \caption{Flow field for mesh ``$48$'' with automatically created matrix.}
    \label{fig:compadmag48AD}
  \end{figure}

\begin{figure}[h!]
    \centering
    % 9.78 x 1.13, width=15.281cm, height=1.7656
    \begin{minipage}[b]{15.281cm}
      \includegraphics[width=15.281cm]{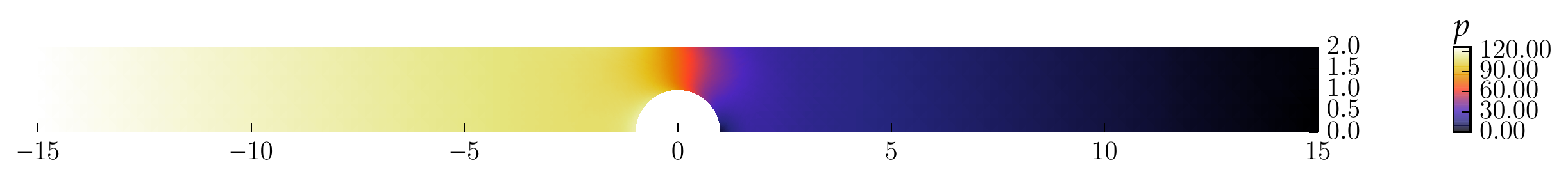}
    \end{minipage}
    \caption{Pressure field for mesh ``$48$'' with automatically created matrix.}
    \label{fig:compadp48AD}
  \end{figure}

  Figures~\ref{fig:compadmag48NOAD}~and~\ref{fig:compadp48NOAD} show the results
  of a simulation with $Wi = 0.05$
  for the rougher mesh
  with a manually derived and implemented matrix.
  A comparison with the results that were obtained using automatic differentiation,
  which are shown in Figures~\ref{fig:compadmag48AD}~and~\ref{fig:compadp48AD},
  does not yield any visible differences.
  This supports the claim that the methods converge to the same solution.

  Some performance issues that arise from the automatic creation of the matrix
  were mentioned in Section~\ref{sec:s3jacfe}, along with one possibility
  for mediating some of these effects. Table~\ref{tab:runtimesmat} shows how
  much time the calculation of the left- and right-hand sides takes up
  for some different mesh resolutions, when the matrix is
  created by different means.
  In agreement with the predictions that were
  made in Section~\ref{sec:s3jacfe}, the execution of the matrix calculation code created
  by automatic differentiation takes much longer than that of the manually created
  code. The ratio between both execution times varies considerably for the different
  simulations.
  When the execution time of the optimized AD code---with separate input vectors
  for all degrees of freedom, as described in Section~\ref{sec:s3jacfe}---is
  compared with that of the unoptimized AD code, the difference between
  the simulations is much less pronounced. The optimized code reliably outperforms
  the unoptimized code, with a performance gain of around $46\%$.

  \begin{table}[h!]
  \centering
  \begin{tabular}{lrrrrr}
  \toprule
  \textbf{Method} & \multicolumn{1}{c}{\textbf{Mesh 48}} & \multicolumn{1}{c}{\textbf{Mesh 96}}
    & \multicolumn{1}{c}{\textbf{Mesh 192}} & \multicolumn{1}{c}{\textbf{Mesh 384}} & \multicolumn{1}{c}{\textbf{Mesh 768}} \\ \midrule
  manual derivative                &   0.98 s &   1.69 s &   3.37 s &   6.63 s &  13.05 s \\
  unoptimized automatic derivative &  13.27 s &  26.48 s &  52.76 s & 105.73 s & 210.28 s \\
  separation by degree of freedom  &   7.04 s &  14.33 s &  28.36 s &  57.24 s & 113.62 s \\
  \midrule
  \textit{number of processors} & \textit{16} & \textit{32} & \textit{64} & \textit{128} & \textit{256} \\
  \bottomrule
  \end{tabular}
  \caption{Time taken up by matrix calculations for different methods (rows) and meshes (columns) with $Wi = 0.05$.}
  \label{tab:runtimesmat}
  \end{table}

  These times can be put into perspective when they are compared to the
  total execution times of the simulation. Table~\ref{tab:runtimesper} shows the percentage of the total
  execution time which is taken up by the left- and right-hand side accumulations in
  the different cases.
  Again, these numbers differ enormously when the different simulations are compared.
  With every increase in the number of processors, the percentage decreases. This suggests
  that the matrix accumulation code scales much better than the linear solver code.
  The performance losses incurred by the method may therefore become negligible for
  some realistic problems.
  However, the numbers in the first two columns show that the influence can become
  very large in certain cases. This confirms the assumption that the optimization of the automatic differentiation
  method remains an important and necessary task.

  \begin{table}[h!]
  \centering
  \begin{tabular}{lrrrrr}
  \toprule
  \textbf{Method} & \multicolumn{1}{c}{\textbf{Mesh 48}} & \multicolumn{1}{c}{\textbf{Mesh 96}}
    & \multicolumn{1}{c}{\textbf{Mesh 192}} & \multicolumn{1}{c}{\textbf{Mesh 384}} & \multicolumn{1}{c}{\textbf{Mesh 768}} \\ \midrule
  manual derivative                & 4.5 \% & 3.5 \% & 1.9 \% & 0.8 \% & 0.3 \% \\
  unoptimized automatic derivative & 40.5 \% & 36.8 \% & 21.0 \% & 10.5 \% & 4.5 \% \\
  separation by degree of freedom  & 30.2 \% & 25.7 \% & 13.5 \% & 6.4 \% & 2.5 \% \\
  \bottomrule
  \end{tabular}
  \caption{Time taken up by matrix calculations for different methods and meshes, given as percentage of total execution time.}
  \label{tab:runtimesper}
  \end{table}

  We want to look at the efficiency of the differentiated code more closely. The
  evaluation of a function resulting from forward mode automatic differentiation
  always includes the evaluation of the original function. Therefore, we can
  assume the computational cost of such an evaluation to be $1 + \alpha$ times
  greater than the cost of the original function evaluation, with positive
  $\alpha$. If the same function is used as the basis for vector forward mode,
  with $p$ directional derivatives being evaluated simultaneously, we can assume
  the cost ratio to be $1 + \alpha \cdot p$.

  In \cite[\S 4.5]{Griewank2008}, it is shown that the upper bound for $\alpha$
  is $1.5$ in this case. For comparison, a finite-difference approach would
  yield $\alpha = 1$.
  For our example, when comparing actual execution times, we could achieve values
  of $\alpha$ between $2.18$ and $2.27$ for the unoptimized approach, and values
  between $1.14$ and $1.23$ for the optimized approach, with separated degrees
  of freedom.

  The values for the unoptimized approach may seem unrealistic at first, when
  compared to the upper bound given in \cite{Griewank2008}. However, we would like
  to point out that these execution times depend greatly on issues such as compiler
  optimization, cache handling, etc., and that the loss in efficiency may also
  partially result from the increased memory requirements of the differentiated
  function.

  One interesting result is that even our optimized approach is still less
  efficient than a one-sided finite-difference approach. Nonetheless, one should
  note that this slight increase in execution time comes with the bonus of a
  much more accurate derivative, and removes the necessity of having to choose
  a suitable step size.

\section{Conclusion}

The objective of this work was to simplify the implementation of partial differential
equations into finite element codes. This was attempted by using automatic differentiation
for the creation of the code for the calculation of the Jacobian matrix.
As the results of the test case in the previous section prove, a useful method
for this purpose could indeed be developed.
However, during the development of this method, several problems surfaced
that could not be solved in a universal manner. In these cases, the solutions
were specific to the program code and equations at hand, and would need to be
adjusted to be useful in different contexts.

One of these problems lies in TAPENADE's incomplete support of the Fortran programming
language.
Even though we used a source code that works perfectly with the popular Fortran compilers,
many changes were necessary in order to make it work with TAPENADE. Some of these changes
could be automated, but it is to be expected that other source codes will bring to light
further problems of this nature.
The complexity of the changes that would need to be automated to support arbitrary Fortran codes
is in fact of such magnitude that only a program equipped with a full-featured Fortran parser
could reliably do the job.

TAPENADE incorporates a feature that is useful for the automatic creation of full Jacobian
matrices instead of single directional derivatives. This is its vector mode for the
simultaneous calculation of arbitrary numbers of directional derivatives of the same
function.
This is a step in a good direction and paves the way for some important optimizations.
However, the resulting code can generally be expected to be less efficient than a
well thought-out manual implementation.
The closed-source nature of TAPENADE means that any further optimization techniques either
have to be applied manually or require a separate language parser. Exploiting sparsity
in input vectors or other properties of the code structure
therefore involves a significant amount of additional work.

One important finding of the present work is the fact that the special structure of
the finite element method allows for an easy automatic creation of a huge Jacobian
matrix with full exploitation of the matrix sparsity by combining the typical finite
element assembly process with automatic differentiation on the element level
(see Section~\ref{sec:s3implprops}).
The final step on the path towards a method that will be ready for production use will now lie in
the combination of this code design with more efficient methods for the automatic creation
of the dense element-level matrix (see Section~\ref{sec:s3jacfe} for some approaches leading in
this direction).
As shown in Section~\ref{sec:s35res}, we already achieved significant efficiency
improvements with a first optimization attempt.

The developments described
in this document had the side effect of prompting the discovery and resolution of a number of
incompatibilities between TAPENADE and current Fortran standards (see Section~\ref{sec:s25incomp}).
The knowledge that was gained this way will aid considerably in the implementation of further projects involving
TAPENADE, that go beyond the calculation of Jacobian matrices.

\section{Acknowledgements}

The authors gratefully acknowledge the support of DFG under the collaborative research project SFB 1120
(subproject B2)
and DFG grant "Computation of Die Swell Behind a Complex Profile Extrusion Die Using a Stabilized Finite
Element Method for Various Thermoplastic Polymers".

%\section*{References}

\bibliography{mybibfile}

\end{document}